\documentclass[a4paper, 12pt] {article}

\font\smallit=cmti10
\begin{document}

\begin{center}
{\bf THE BOUNIAKOWSKY CONJECTURE AND THE DENSITY OF POLYNOMIAL ROOTS TO 
PRIME MODULI}
\vskip 20pt
{\bf Timothy Foo}\\
{\smallit Division of Mathematical Sciences, Nanyang Technological 
University, Singapore}\\
{\tt S080074@ntu.edu.sg}
\end{center}
\vskip 30pt

\noindent
\begin {abstract} We establish a result linking the Bouniakowsky 
conjecture 
and the density of polynomial roots to prime moduli.
\end{abstract}
\vskip 5pt



\section*{\normalsize Introduction}
In this paper, we study roots of irreducible polynomials to prime moduli. 
We think of $\mathbf{Z}/p\mathbf{Z}$ as the set ${0,1,2,...,p-1}$ and 
hence we think of the root of our polynomial as a number in that set. When 
the root $z$ is divided by $p$, we naturally have a number in $(0,1)$. If 
we fix a polynomial $f(x)$ of degree $n \geq 2$ which is irreducible in 
$\mathbf{Z}[x]$, we can 
consider the set $$A_f = \cup_p \{\frac{z}{p}: f(z) \equiv 0 \mbox{ mod 
}p, 
1 
\leq z 
\leq p-1\}$$ The aim of this paper is to prove that if a certain 
conjecture called the Bouniakowsky conjecture is true, then the set $A_f$ 
is dense in (0,1). We stress that our result is conditional. Results that 
are not dependent on open conjectures have been proven about roots of 
polynomials to various moduli. Hooley [H] proved that the roots of an 
irreducible polynomial, considered over the ring $\mathbf{Z}/n\mathbf{Z}$, 
$n$ not necessarily prime, when suitably normalized by dividing by $n$ and 
considered over all $n$, are in fact equidistributed in $(0,1)$. Duke, 
Friedlander and Iwaniec[DFI] proved equidistribution for quadratic 
polynomials of negative discriminant, to prime moduli. Toth[T] proved 
equidistribution for quadratic polynomials of positive discriminant, to 
prime moduli. We now state the main theorem of our paper.

\vskip 5pt
\noindent
{\bf Theorem} If the Bouniakowsly conjecture is true, the set $A_f = 
\cup_p \{\frac{z}{p}: f(z) \equiv 0 \mbox{ 
mod } p, 1 \leq z \leq p-1\}$ is dense in (0,1). 

\vskip 20pt
\section*{\normalsize The Bouniakowsky conjecture} 
We now discuss the Bouniakowsky conjecture to give some background. 

\vskip 5pt
\noindent
{\bf Conjecture} Bouniakowsky Conjecture: Let $f(x)$ be a polynomial that 
is irreducible in $\mathbf{Z}[x]$. Let $r_f = \gcd(\{f(x): x \in 
\mathbf{Z}\})$. Then $\frac{f(x)}{r_f}$ is prime infinitely often. 

\vskip 5pt
\noindent
It is easy to construct polynomials which are always divisible by a given 
prime $q$. We know by Fermat's little theorem that the prime $q$ always 
divides $x^q-x$. Therefore, all we have to do is choose a value $k$ so 
that $x^q - x + qk$ is irreducible in $\mathbf{Z}[x]$. It then follows 
that $q$ divides all the values of this polynomial.  

\vskip 20pt
\section*{\normalsize The Result}
We first begin by considering a subset of $(0,1)$ which we will prove 
to be dense. We are then going to use this set to help prove the density 
of $A_f$. Here, we let $n$ be the degree of $f$ and $c$ be the leading 
coefficient of $f$.

\vskip 5pt
\noindent
Let $B_f = \{\frac{a}{b}: 1 \leq a < b, b \mbox{ odd prime, } 
(cr_f,b)=1, acx^{n-1} = -r_f \mbox{ mod b has a solution} \}$

\vskip 5pt
\noindent
{\bf Lemma 1} $B_f$ is dense in $(0,1)$.

\vskip 5pt
\noindent
{\it Proof}
 
\noindent
Case 1: $n$ is even. Consider the map $x \rightarrow x^{n-1}$ on 
$(\mathbf{Z}/b\mathbf{Z})^{*}$. This map is injective and surjective if 
$(n-1,b-1) = 1$. For such $b$, we can in fact solve $acx^{n-1} \equiv 
-r_f$ mod $b$ for all $a \in (\mathbf{Z}/b\mathbf{Z})^{*}$. Since $b$ is 
prime, we can pick $b$ larger than $cr_f$ to ensure $(b,cr_f)=1$. We can 
also pick infinitely many such $b$ with $(n-1,b-1)=1$. It thus follows 
that $B_f$ is dense in this case.

\vskip 5pt
\noindent
Case 2: $n$ is odd. Since $n-1$ is even, let $n-1=2^eh$, $h$ odd. The map 
$x \rightarrow x^{n-1}$ on $(\mathbf{Z}/b\mathbf{Z})^{*}$ is therefore a 
composition of the maps $x \rightarrow x^2$ applied $e$ times and $x 
\rightarrow x^h$. $x \rightarrow x^h$ is a permutation of 
$(\mathbf{Z}/b\mathbf{Z})^{*}$ if $(b-1,h)=1$. Also, if $b \equiv 3$ mod 
$4$, $x \rightarrow x^2$ is a permutation of the squares in 
$(\mathbf{Z}/b\mathbf{Z})^{*}$, so by choosing $b \equiv 3$ mod $4$ and 
$(b-1,h) = 1$, we can ensure that the image of $x \rightarrow x^{n-1}$ is 
the squares. We also want $(b,cr_f)=1$. We have infinitely many primes $b$ 
satisfying these conditions, and for such $b$, the numerator of the 
fractions $\frac{a}{b}$ ranges over either only the squares or only the 
nonsquares in $(\mathbf{Z}/b\mathbf{Z})^{*}$. By a result of Brauer, the 
maximum number of consecutive squares or nonsquares in 
$(\mathbf{Z}/b\mathbf{Z})^{*}$ is less than $b^{0.5}$ when $b \equiv 3$ 
mod 
$4$. [B]
This ensures that $B_f$ is dense in this case. 

\vskip 5pt
\noindent
We will now show how $\frac{z}{p}$ is related to the values in $B_f$. To 
do this, first consider the original polynomial $f$. From $f = 
\sum_{i}c_ix^i$, we can construct a polynomial $g(x,y) = 
\sum_{i}c_ix^iy^{n-i}$. Now for any prime $b$ with $(b,cr_f)=1$ we have a 
polynomial in one variable $g(bw+t,b)$ where $w$ is the variable and 
$t\in (\mathbf{Z}/b\mathbf{Z})^{*}$. Since we can vary $b$ and $t$, we 
have many such polynomials associated to $f$. We will show that the gcd 
of the values of all these polynomials is also $r_f$ and that they are 
also irreducible in $\mathbf{Z}[w]$. It is these 
polynomials that we apply the Bouniakowsky conjecture to. If the 
Bouniakowsky conjecture is true, then there are infinitely many primes 
$p$ with $r_fp = g(bw+t,b)$ as $w \rightarrow \infty$. 
Furthermore, for 
these primes $p$, we can construct a root $z$ of $f$ mod $p$ such that 
$\frac{z}{p}$ is ``close" to $\frac{a}{b}$ where $a$ is chosen so that 
$\frac{ap+bw+t}{b}$ is an integer and $\frac{a}{b} \in (0,1)$. This is the 
same as choosing $1 \leq a < b$ and $a$ such that $act^{n-1} \equiv -r_f$ 
mod $b$. We thus see the relation to the set $B_f$. We then let 
$z=\frac{ap+bw+t}{b}$ and show 
that $z$ is a root of $f$ mod $p$.

\vskip 5pt
\noindent
{\bf Lemma 2} The polynomial $g(bw+t,b)$, where $w$ is the variable, $b$ 
is prime, $(b,cr_f)=1$, $1 \leq t < b$ is irreducible in $\mathbf{Z}[w]$.

\vskip 5pt
\noindent
{\it Proof} $g(bw+t,b)$ is related in a simple way to the original 
polynomial $f$. $g(bw+t,b) = \sum_ic_i(bw+t)^ib^{n-i} = 
b^n\sum_ic_i(w+\frac{t}{b})^i = b^ng(w+\frac{t}{b},1) = 
b^nf(w+\frac{t}{b})$. Since a polynomial is irreducible in $\mathbf{Z}[x]$ 
if and only if it is irreducible in $\mathbf{Q}[x]$, the lemma follows.

\vskip 5pt
\noindent
{\bf Lemma 3} Let $b$ be prime, $(b,cr_f)=1$, and $1 \leq t < b$. Then 
$\gcd(\{g(bw+t,b): w \in \mathbf{Z} \})=r_f$

\vskip 5pt
\noindent
{\it Proof} Let $r=r_f$. Since $f$ has integer coefficients, we can think 
of $f$ as a polynomial in $(\mathbf{Z}/r\mathbf{Z})[x]$. But since $r$ 
divides all the values of $f$, it follows that $f(x) = 0$ in 
$(\mathbf{Z}/r\mathbf{Z})[x]$. We showed in the proof 
of 
Lemma 2 that $g(bw+t,b) = b^nf(w+\frac{t}{b})$ in $\mathbf{Q}[x]$. Since 
$(b,r_f)=1$, $b$ has an inverse mod $r$ and hence the rational number 
$\frac{t}{b}$ can be thought of as an element in $\mathbf{Z}/r\mathbf{Z}$. 
Hence $g(bw+t,b) = b^nf(w+\frac{t}{b}) = 0$ 
in $(\mathbf{Z}/r\mathbf{Z})[x]$. Therefore, for each such $b$ and $t$, we 
have that $r$ divides $\gcd(\{g(bw+t,b): w \in \mathbf{Z}\})$.
Conversely, let $r_{b,t} = \gcd(\{g(bw+t,b): w \in \mathbf{Z}\})$. We have 
$g(bw+t,b) = 0$ in $(\mathbf{Z}/r_{b,t}\mathbf{Z})[w]$. But $f(w) = 
(b^n)^{-1}g(b(w-\frac{t}{b})+t,b)$, so $f(w) = 0$ in 
$(\mathbf{Z}/r_{b,t}\mathbf{Z})[w].$ Therefore $r_{b,t}$ divides $r$ for 
each such $b$ and $t$. 
It follows that the polynomials $g(bw+t,b)$ 
have the same gcd as $f$.

\vskip 5pt
\noindent
{\bf Lemma 4} If $a$ is chosen such that $z = \frac{ap+bw+t}{b}$ is an 
integer, then $z$ is a root of the polynomial $f$ mod 
$p$.

\vskip 5pt
\noindent
{\it Proof} $b^nf(z) = b^nf(\frac{ap+bw+t}{b}) = 
b^n\sum_ic_i(\frac{ap+bw+t}{b})^i = \sum_ic_i(ap+bw+t)^ib^{n-i} \equiv 
\sum_ic_i(bw+t)^ib^{n-i} = g(bw+t,b) = r_fp \equiv 0$ mod $p$. Since 
$(b,p)=1$, the lemma is proven.

\vskip 5pt
\noindent
Having proven these lemmas, we know that $\frac{z}{p}$ is close to 
$\frac{a}{b}$. 
Assuming the 
Bouniakowsky conjecture, we can let $w \rightarrow \infty$ and obtain 
infinitely many primes $p$ and a 
root $z$ for each prime. As $w \rightarrow \infty$, $\frac{z}{p}$ is 
arbitrarily close to $\frac{a}{b}$, since $n \geq 2$. Since we showed in 
Lemma 1 that $B_f$ 
is dense in $(0,1)$, the theorem is now proved. 

\vskip 5pt
\noindent
\section*{\normalsize Acknowledgements}
Thanks to Professor Zhengyu Mao for suggesting the problem to me and for 
many helpful discussions.

\vskip 5pt
\noindent
\section*{\normalsize References}
[B] A. Brauer, \"{U}ber die Verteilung der Potenzreste, Math Z. 35 (1932) 
39-50.

\noindent
[DFI] W. Duke, J.B. Friedlander, H. Iwaniec, Equidistribution of roots of 
a quadrtic congruence to prime moduli, Annals of Mathematics (2) 141 
(1995), no. 2, 423-441.

\noindent
[H] C. Hooley, On the distribution of the roots of a polynomial 
congruences, Mathematika 11 (1964) 39-49

\noindent
[T] A. Toth, Roots of quadratic congruences, IMRN 2000, no. 14, 719-739
  
\end{document}